\documentclass[11pt]{article}
\usepackage{amsmath}
\usepackage{amssymb}
\usepackage{latexsym}
\begin{document}
\title{Hom-Akivis algebras}
\author{\sc {A. Nourou Issa } \\ 
\sc {D\'epartement de Math\'ematiques, Universit\'e d'Abomey-Calavi}, \\
01 BP 4521 \sc {Cotonou 01, B\'enin} \\ 
{\it E-mail}: woraniss@yahoo.fr}
\date{}
\maketitle
\begin{abstract}
Hom-Akivis algebras are introduced. The commutator-Hom-associator algebra of a non-Hom-associative algebra (i.e. a Hom-nonassociative algebra) is a Hom-Akivis algebra. It is shown that non-Hom-associative algebras can be obtained from nonassociative algebras by twisting along algebra automorphisms while Hom-Akivis algebras can be obtained from Akivis algebras by twisting along algebra endomorphisms. It is pointed out that a Hom-Akivis algebra associated to a Hom-alternative algebra is a Hom-Malcev algebra.
\\
{\it Classification}: 17A30, 17D10, 17D99
\end{abstract}
{\bf 1. Introduction}
\\
\par
The theory of Hom-algebras originated from Hom-Lie algebras introduced by J.T. Hartwig, D. Larsson, and S.D. Silvestrov in [9] in the study of quasi-deformations of Lie algebras of vector fields, including q-deformations of Witt algebras and Virasoro algebras. The connection between the theory of Hom-algebras and deformation theory and other trends in mathematics attracted attention of researchers (see, e.g., [5], [6], [7], [8], [13], [14], [20]). Generalizing the relation between Lie algebras and associative algebras, the notion of a Hom-associative algebra is introduced by A. Makhlouf and S.D. Silvestrov in [12], where it is shown that the commutator algebra (with the twisting map) of a Hom-associative algebra is a Hom-Lie algebra. By twisting defining identities, other Hom-type algebras such as Hom-alternative algebras, Hom-Jordan algebras [14], Hom-Novikov algebras [20], or Hom-Malcev algebras [21] are introduced and discussed.
\par
As for Hom-alternative algebras or Hom-Novikov algebras, we consider in this paper a twisted version of the Akivis identity which defines the so-called Akivis algebras. We call ``Hom-Akivis algebra'' this twisted Akivis algebra. It is known [3] that the commutator-associator algebra of a nonassociative algebra is an Akivis algebra. This led us to consider ``non-Hom-associative algebras'' i.e. Hom-nonassociative algebras or nonassociative Hom-algebras ([12], [18]) and we point out that the commutator-Hom-associator algebra of a non-Hom-associative algebra has a Hom-Akivis structure. In this setting, Akivis algebras are special cases of Hom-Akivis algebras in which the twisting map is the identity map. Also the class of Hom-Akivis algebras contains the one of Hom-Lie algebras in the same way as the class of Akivis algebras contains the one of Lie algebras.
\par
Akivis algebras were introduced by M. A. Akivis ([1], [2], [3]) as a tool in the study of some aspects of web geometry and its connection with loop theory. These algebras were originally called ``$W$-algebras'' [3]. Later, Hofmann and Strambach [10] introduced the term ``Akivis algebras'' for such algebraic objects.
\par
The rest of the present paper is organized as follows. In Section 2 we recall basic definitions and useful results about Akivis algebras, Hom-Lie algebras and Hom-associative algebras. In Section 3 we consider non-Hom-associatve algebras (one observes the counterpart of the generalization of associative algebras by nonassociative ones). We point out that multiplicative non-Hom-associative algebras are constructed from nonassociative algebras via automorphisms (Theorem 3.3) while, in the case of Hom-associative and associative algebras, one needs just endomorphisms [19]. In Section 4 Hom-Akivis algebras are considered. Two methods of producing Hom-Akivis algebras are provided starting with either non-Hom-associative algebras (Theorem 4.2) or classical Akivis algebras along with twisting maps (Corollary 4.6). Hom-Akivis algebras are shown to be closed under twisting by self-morphisms (Theorem 4.4). In Section 5, Hom-Akivis algebras associated to Hom-alternative algebras are shown to be Hom-Malcev algebras (these later algebraic objects are recently introduced by D. Yau [21]). This could be seen as a generalization of the Malcev construction of Moufang-Lie algebras (i.e. Malcev algebras) from alternative algebras [15]. 
\par
Throughout this paper, all vector spaces and algebras are meant over a ground field $\mathbb{K}$ of characteristic 0. \\
\\
{\bf 2. Preliminaries}
\\
\par
We recall useful definitions and results that are for further use. We begin with Akivis algebras.
\par
An {\it Akivis algebra} ($L, [-,-], \langle -, -, - \rangle$) is a vector space $L$ together with a bilinear skew-symmetric binary operation $(x, y) \mapsto [x,y]$ and a trilinear ternary operation $(x, y, z) \mapsto \langle x, y, z \rangle$ that are linked by the identity \\
\par
$\sigma [[x,y],z] = \sigma  \langle x, y, z \rangle - \sigma  \langle y, x, z \rangle$, \hfill (2.1) \\
\\
where here, and in the sequel, $\sigma $ denotes the sum over cyclic permutation of $x, y, z$.
\par
The relation (2.1) is called the {\it Akivis identity}.
\par
In loop theory, roughly, Akivis algebras are for local smooth loops what are Lie algebras for local Lie groups. However Akivis algebras originated from web geometry ([1], [2]; see also [4] for a survey of the subject).
\par
It is well known that the commutator algebra of an associative algebra is a Lie algebra. M. A. Akivis [3] generalized this construction to nonassociative algebras. Indeed, he pointed out the following \\
\par
{\bf Theorem 2.1.} [3] {\it The commutator-associator algebra of a nonassociative algebra is an Akivis algebra}. \\
\par
The Akivis algebra constructed by Theorem 2.1 is said {\it associated} (with a given nonassociative algebra) [10].
\par
In Section 4 (Theorem 4.2) we give the Hom-counterpart of Theorem 2.1. Beforehand we recall some facts about Hom-algebras.
\par
A {\it Hom-module} [19] is a pair ($ V, {\alpha}_{V} $), where $V$ is a vector space and $ {\alpha}_{V} : V \mapsto V $ a linear map.
\par
A {\it Hom-associative algebra} [12] is a triple ($ V, {\mu}_{V}, {\alpha}_{V} $) in which ($ V, {\alpha}_{V} $) is a Hom-module and $ {\mu}_{V} : V \times V \mapsto V $ is a bilinear operation on $V$ such that \\
\par
$ {\mu}_{V} ({\mu}_{V}(x,y), {\alpha}_{V}(z)) =  {\mu}_{V}( {\alpha}_{V}(x), {\mu}_{V} (y,z))$, \hfill (2.2) \\
\\
for $ x,y,z \in V$. 
\par
The relation (2.2) is called the {\it Hom-associativity} for ($ V, {\mu}_{V}, {\alpha}_{V} $). If ${\alpha}_{V} = {id}_{V}$, then (2.2) is just the associativity.
\par
Using the abbreviation $xy$ for ${\mu}_{V}(x,y)$, the Hom-associativity (2.2) reads \\
\par
$(xy) {\alpha}_{V}(z) = {\alpha}_{V}(x)(yz)$. \\
\par
The Hom-associative algebra ($ V, {\mu}_{V}, {\alpha}_{V} $) is said {\it multiplicative} if ${\alpha}_{V}$ is an endomorphism (i.e. a self-morphism) of ($ V, {\mu}_{V}$). Hom-associative algebras are closely related to Hom-Lie algebras.
\par
A {\it Hom-Lie algebra} is a triple ($ V, [-,-], {\alpha}_{V} $) in which ($ V, {\alpha}_{V} $) is a Hom-module and $ [-,-] : V \times V \mapsto V $ is a bilinear skew-symmetric operation on $V$ such that \\
\par
$\sigma [[x,y], {\alpha}_{V}(z)] = 0 $, \hfill (2.3) \\
\\
for $ x,y,z \in V$.
\par
The relation (2.3) is called the {\it Hom-Jacobi identity}. If, moreover, ${\alpha}_{V}$ is an endomorphism of ($ V, [-,-] $), then ($ V, [-,-], {\alpha}_{V} $) is said {\it multiplicative}. Examples of Hom-Lie algebras could be found in [9], [12], [19].
\par
In the Hom-Lie setting, the Hom-associative algebras play the role of associative algebras in the Lie setting in this sense that the commutator algebra of a Hom-associative algebra is a Hom-Lie algebra [12]. In [19] it is shown how arbitrary associative (resp. Lie) algebras give rise to Hom-associative (resp. Hom-Lie) algebras via endomorphisms. These constructions are considered here in the Hom-Akivis setting. \\
\\
{\bf 3. Non-Hom-associative algebras. Examples and construction} \\
\par
In [12], [18] the notion of a Hom-associative algebra is extended to the one of a Hom-nonassociative algebra (or nonassociative Hom-algebra, or just Hom-algebra), i.e. a Hom-type algebra in which the Hom-associativity (as defined by (2.2)) does not necessarily hold. In this section and in the rest of the paper, such a Hom-type algebra will be called ``non-Hom-associative'' (although this terminology seems to be somewhat cumbersome) in order to stress the Hom-counterpart of the generalization of associative algebras by the nonassociative ones. We provide some examples and show that a  specific twisting of a nonassociative algebra gives rise to a non-Hom-associative algebra. \\
\par
{\bf Definition 3.1.} ([12], [18]). A {\it multiplicative non-Hom-associative} (i.e. not necessarily Hom-associative) {\it algebra} is a triple ($A, \mu , \alpha $) such that \\
(i) ($A, \alpha $) is a Hom-module; \\
(ii) $ \mu: A \times A \mapsto A $ is a bilinear operation on $A$; \\
(iii) $\alpha$ is an endomorphism of ($A, \mu $) (multiplicativity). \\
\par
If $\alpha$ is the identity map in Definition 3.1, then ($A, \mu , \alpha $) reduces to a nonassociative algebra ($A, \mu $). 
\par
The following example of a non-Hom-associative algebra is derived from an example in [14]. \\
\par
{\bf Example 3.2.} Let $\{ u, v, w \}$ be a basis of a three-dimensional vector space $V$. Define on $V$ the operation $\mu$ and the linear map $\alpha$ as follows:
\begin{displaymath}
\begin{array}{clcr}
&\mu (u,u) = au,               &\mu (v,v) =  av  & \\
&\mu (u,v) = \mu (v,u) = av,   &\mu (v,w) =  bw  & \\
&\mu (u,w) = \mu (w,u) = bw,   &\mu (w,v) =  &\mu (w,w) = 0, 
\end{array}
\end{displaymath}
and
\par
$\alpha (u) = av$, $\alpha (v) = aw$, $\alpha (w) = bu$, \\
where $a,b \in \mathbb{K}$ and $a \neq 0$, $b \neq 0$.
\par
Now we have
\par
$\mu (\mu (u,v), \alpha (w)) = \mu (av, bu) = a^{2}bv$ \\
and
\par
$\mu ( \alpha (u), \mu (v,w)) = \mu (av, bw) = ab^{2}w$ \\
so $ \mu (\mu (u,v), \alpha (w)) \neq \mu ( \alpha (u), \mu (v,w))$, i.e. the Hom-associativity fails for the triple ($V, \mu , \alpha $), and thus ($V, \mu , \alpha $) is non-Hom-associative. Observe that if choose
$a \neq b$ (with $a \neq 0$, $b \neq 0$), then ($V, \mu $) is nonassociative since $ \mu (\mu (v,v), w) \neq \mu (v, \mu (v,w))$. \\
\par
Other examples follow. \\
\par
{\bf Example 3.3.} The Lie algebra ${\mathfrak sl}(2, {\mathbb C})$ has a basis $\{ u, v, w \}$ with multiplication: \\
\par
$[u,v] = -2u$,  $[u,w] = v$, $[v,w] = -2w$. \\
\par
Define $\alpha : {\mathfrak sl}(2, {\mathbb C}) \rightarrow {\mathfrak sl}(2, {\mathbb C})$ by setting: \\
\par
$\alpha (u) = w$, $\alpha (v) = -v$, $\alpha (w) = u$. \\
\par
Then $\alpha$ is a self-morphism of $({\mathfrak sl}(2, {\mathbb C}), [-,-])$. Next $[\alpha (u), [w,w]] = 0$ \\ while $[[u, w],\alpha (w)] = 2u$ so that $({\mathfrak sl}(2, {\mathbb C}), [-,-], \alpha)$ is non-Hom-associative. \\
\par
{\bf Example 3.4.} There is a five-dimensional nonassociative (flexible) algebra $(A, \cdot )$ with basis $\{ e_{1}, ..., e_{5} \}$ and multiplication: \\
\par
$e_{1} \cdot e_{2} = e_{5} + \frac{1}{2} e_{4}$, \;\;\; $e_{1} \cdot e_{4} = \frac{1}{2} e_{1} = - e_{4} \cdot e_{1}$,
\par
$e_{2} \cdot e_{1} = e_{5} - \frac{1}{2} e_{4}$, \;\;\; $e_{2} \cdot e_{4} = - \frac{1}{2} e_{2} = - e_{4} \cdot e_{2}$,
\par
$e_{3} \cdot e_{4} = \frac{1}{2} e_{3} = - e_{4} \cdot e_{3}$,
\par
$e_{4} \cdot e_{4} = - e_{5}$, \\
\\
and all other products are $0$ (see [16], p. 29, Example 1.5). Define $\alpha : A \rightarrow A$ by: \\
\par
$\alpha (e_{1}) = e_{2}$, \;\; $\alpha (e_{2}) = e_{1}$, \;\; $\alpha (e_{3}) = 0$, \;\; $\alpha (e_{4}) = - e_{4}$, \;\; $\alpha (e_{5}) = e_{5}$. \\
\\
Then $\alpha$ is a self-morphism of $(A, \cdot )$. Moreover $(A, \cdot , \alpha)$ is non-Hom-associative since $ \alpha (e_{3}) \cdot e_{4}e_{4} = 0$ while $e_{3}e_{4} \cdot \alpha (e_{4}) = - \frac{1}{4} e_{3}$. However, by Theorem 4.4 in [21], $(A, \cdot , \alpha)$ is Hom-flexible. \\
\par
{\bf Example 3.5.} The commutator algebra $A^{-}$ of the algebra $A$ of Example 3.4 is defined by: \\
\par
$e_{1} \star e_{2} = e_{4} = - e_{2} \star e_{1}$, \;\; $e_{1} \star e_{4} = e_{1} = - e_{4} \star e_{1}$,
\par
$e_{2} \star e_{4} = - e_{2} = - e_{4} \star e_{2}$, 
\par
$e_{3} \star e_{4} = e_{3} = - e_{4} \star e_{3}$, \\
\\
and all other products are $0$. Then $(A^{-}, \star )$ is a Malcev algebra ([16], p. 29, Example 1.5) and, as observed in [16], $(A^{-}, \star )$ is isomorphic to the unique five-dimensional non-Lie solvable Malcev algebra found by E. N. Kuzmin [11]. Next, define  the map $\alpha : A^{-} \rightarrow A^{-}$ as in Example 3.4 above. Then $\alpha$ is a self-morphism of $(A^{-}, \star )$ and $(A^{-}, \star , \alpha)$ is non-Hom-associative since $ \alpha (e_{3}) \star e_{4}e_{4} = 0$ and $e_{3}e_{4} \star \alpha (e_{4}) = - e_{3}$. \\
\par
The examples above show that the nonassociativity of ($V, \mu $) is not enough for the failure of Hom-associativity in ($V, \mu , \alpha $) and thus the choice of the twisting map is of prime significance in the definition of Hom-type algebras. In other words, there is a freedom on how to twist. This is observed and investigated in [6]. In the present paper, we consider only twisting maps which ensures the failure of the Hom-associativity as defined by (2.2).
\par
The following theorem could be seen as the nonassociative counterpart of the construction of Hom-associative algebras from associative algebras given by D. Yau in [19]. But here an additional condition on the twisting map is needed. \\
\par
{\bf Theorem 3.6.} {\it Let $(A, \mu )$ be a nonassociative algebra and let $\alpha$ be an automorphism of $(A, \mu )$. Then $(A, {\mu}_{\alpha} = \alpha \circ \mu , \alpha )$ is a multiplicative non-Hom-associative algebra.
\par
Moreover, let $(\tilde{A}, \tilde{\mu} )$ be another nonassociative algebra and $\tilde{\alpha}$ an automorphism of $(\tilde{A}, \tilde{\mu} )$. If $f: A \rightarrow \tilde{A}$ is an algebra morphism that satisfies $f \circ \alpha = \tilde{\alpha} \circ f$, then $ f: (A, {\mu}_{\alpha}, \alpha ) \rightarrow (\tilde{A}, \tilde{\mu}_{\tilde{\alpha}} = \tilde{\alpha} \circ \tilde{\mu} , \tilde{\alpha} )$ is a morphism of multiplicative non-Hom-associative algebras}. \\
\par
{\sc Proof}: Denote $\mu (x,y) := xy$. Then we have \\
${\mu}_{\alpha}({\mu}_{\alpha}(x,y), \alpha (z)) = (\alpha \circ \mu)((\alpha \circ \mu (x,y),\alpha (z)) = \alpha (\alpha (xy)\alpha (z)) \\ = \alpha ((\alpha (x)\alpha (y))\alpha (z)) $ \\
and \\
${\mu}_{\alpha} (\alpha (x), {\mu}_{\alpha}(y,z)) = (\alpha \circ \mu)(\alpha (x), (\alpha \circ \mu)(y,z)) = \alpha ( \alpha (x) \alpha (yz)) \\ = \alpha (\alpha (x) (\alpha (y) \alpha (z)))$. \\
The nonassociativity of ($A, \mu $) and the surjectivity of $\alpha$ imply that \\ $ (\alpha (x)\alpha (y))\alpha (z) \neq \alpha (x) (\alpha (y) \alpha (z))$ so that, by the injectivity of $\alpha$, \\ ${\mu}_{\alpha}({\mu}_{\alpha}(x,y), \alpha (z)) \neq {\mu}_{\alpha} (\alpha (x), {\mu}_{\alpha}(y,z))$ with $x,y,z$ in $A$. Thus the Hom-associativity needs not to hold in $(A, {\mu}_{\alpha}, \alpha )$.
\par
Next, ${\mu}_{\alpha} (\alpha (x), \alpha (y)) = (\alpha \circ \mu)(\alpha (x), \alpha (y)) = \alpha ((\alpha (x)\alpha (y)) = \alpha (\alpha (xy)) = \alpha ({\mu}_{\alpha}(x,y))$ so that $\alpha$ is multiplicative with respect to ${\mu}_{\alpha}$.
\par
The second part of the theorem is proved as follows. \\
$f({\mu}_{\alpha}(x,y)) = f(\alpha (xy)) = (f \circ \alpha) (xy) =(\tilde{\alpha} \circ f)(xy) = \tilde{\alpha} (f(xy))  \\ = \tilde{\alpha}(\tilde{\mu}(f(x),f(y))) = (\tilde{\alpha} \circ \tilde{\mu})(f(x),f(y)) = \tilde{\mu}_{\tilde{\alpha}}(f(x),f(y))$. \hfill $\square$ \\
\par
{\bf Remark.} For a non-Hom-associative algebra ($A, \mu, \alpha $), the automorphism property of $\alpha$ ensures the existence of a nonassociative structure ($A, {\mu}'$) on $A$ which induces the given non-Hom-associative structure ($A, \mu, \alpha $) by Theorem 3.6. In fact, ${\mu}'$ is the untwisted operation on $A$ and ${\mu}' = {\alpha}^{-1} \circ \mu $, since ${\alpha}^{-1}$ is also an algebra automorphism. Such an observation is made in [14] for an alternative algebra ($V, \mu $) and its induced Hom-alternative algebra ($V, {\mu}_{\alpha}, \alpha $) in case when $\alpha$ is an algebra automorphism. The same question is discussed in [7], [8] for associative and Hom-associative algebras. \\
\par
{\bf Example 3.7.} Let $R$ be a unital nonassociative algebra over $\mathbb{K}$ and let $R_{n}$ denote the algebra of $n \times n$ matrices with entries in $R$ and matrix multiplication denoted by $\mu (x,y) = xy$. Then $A := (R_{n}, \mu)$ is also a unital nonassociative algebra. Denote by $N(A)$ the nucleus of $A$ and suppose that there is an invertible element $u \in N(A)$ and $u^{-1} \in N(A)$. Then the map $\alpha (u): A \rightarrow A$ defined by $\alpha (u)(x) = uxu^{-1}$ for $x \in R_{n}$, is an automorphism of $A$.
\par
Define ${\mu}_{\alpha (u)} (x,y) = u(xy)u^{-1}$, for all $x,y \in R_{n}$. Then one checks that $A_{u} = (R_{n} , {\mu}_{\alpha (u)} , \alpha (u))$ is a multiplicative non-Hom-associative algebra and $\alpha (u)$ is an automorphism of $(R_{n} , {\mu}_{\alpha (u)})$. In this way we get a family \{$A_{u} : u \in R_{n} \; \mbox {invertible and} \; u, u^{-1} \in N(A)$\} of multiplicative non-Hom-associative algebras. \\
\par
An example, similar to the example 3.7 above, is given in [19] describing Hom-associative deformations by inner automorphisms.
\par
For a non-Hom-associative algebra, it makes sense in considering the so-called Hom-associators [14] (see also [12]), just as associators are considered in a nonassociative algebra.
\par
Let ($V, \mu, \alpha $) be a non-Hom-associative algebra, where $V$ is a $\mathbb{K}$-linear space, $\mu$ a bilinear operation on $V$ and $\alpha$ a twisting map. For any $x,y,z \in V$, the {\it Hom-associator} is defined by \\
\par
$\mathfrak{as} (x,y,z) = \mu (\mu (x,y), \alpha (z)) - \mu ( \alpha (x) , \mu (y,z))$. \hfill (3.1) \\
\par
Then ($V, \mu, \alpha $) is said:
\par
{\it Hom-flexible}, if $\mathfrak{as} (x,y,x) = 0$ ;
\par
{\it Hom-alternative}, if $\mathfrak{as} (x,y,z)$ is skew-symmetric in $x,y,z$. \\
\\
{\bf 4. Hom-Akivis algebras. Construction} \\
\par
In this Section we give the notion of a (multiplicative) Hom-Akivis algebra that could be seen as a generalization of an Akivis algebra and we point out that such a notion does fit with the one of a non-Hom-associative algebra given in Section 3. In fact we prove the analogue of the Akivis construction (see Theorem 2.1) that the commutator-Hom-associator algebra of a given non-Hom-associative algebra is a Hom-Akivis algebra (Theorem 4.2). Moreover, following [19] for Hom-associative algebras and Hom-Lie algebras, we give a procedure for the construction of Hom-Akivis algebras from Akivis algebras and their algebra endomorphisms. \\
\par
{\bf Definition 4.1.} A {\it Hom-Akivis algebra} is a quadruple \\ ($V, [-,-], [-,-,-], \alpha$), where $V$ is a vector space, $[-,-] : V \times V \rightarrow V$ a skew-symmetric bilinear map, $[-,-,-] : V \times V \times V \rightarrow V$ a trilinear map and $\alpha : V \rightarrow V$ a linear map such that \\
\par
$\sigma [[x,y], \alpha (z)] = \sigma [x,y,z] - \sigma [y,x,z]$, \hfill (4.1) \\
\\
for all $x,y,z$ in $V$.
\par
A Hom-Akivis algebra ($V, [-,-], [-,-,-], \alpha$) is said multiplicative if $\alpha$ is an endomorphism with respect to $[-,-]$ and $[-,-,-]$. \\
\par
In analogy with Lie and Akivis cases, let call (4.1) the {\it Hom-Akivis identity}. \\
\par
{\bf Remark.} (1) If $\alpha = {id}_{V}$, the Hom-Akivis identity (4.1) is the usual Akivis identity (2.1).
\par
(2) The Hom-Akivis identity (4.1) reduces to the Hom-Jacobi identity (2.3), when $[x,y,z] = 0$, for all $x,y,z$ in $V$. \\
\par
The following result shows how one can get Hom-Akivis algebras from non-Hom-associative algebras. \\
\par
{\bf Theorem 4.2.} {\it The commutator-Hom-associator algebra of a multiplicative non-Hom-associative algebra is a multiplicative Hom-Akivis algebra}. \\
\par
{\sc Proof}: Let ($A, \mu , \alpha $) be a multiplicative non-Hom-associative algebra. For any $x,y,z$ in $A$, define the operations
\par
$ [x,y] := \mu (x,y) - \mu (y,x)$ (commutator)
\par
$[x,y,z]_{\alpha} := \mathfrak{as} (x,y,z)$ (Hom-associator; see (3.1)).
\par
For simplicity, set $xy$ for $\mu (x,y)$. Then
\par
$[[x,y], \alpha (z)] = (xy) \alpha (z) - (yx) \alpha (z) - \alpha (z) (xy) + \alpha (z) (yx)$ \\
and
\par
$[x,y,z]_{\alpha} - [y,x,z]_{\alpha} = (xy) \alpha (z) - \alpha (x) (yz ) - (yx) \alpha (z) + \alpha (y)(xz)$. \\
Expanding $\sigma [[x,y], \alpha (z)]$ and $\sigma ([x,y,z]_{\alpha} - [y,x,z]_{\alpha})$ respectively, one gets (4.1) and so ($A, [-,-], [-,-,-]_{\alpha}, \alpha$) is a Hom-Akivis algebra. The multiplicativity of 
($A, [-,-], [-,-,-]_{\alpha}, \alpha$) follows from the one of ($A, \mu , \alpha $). \hfill $\square$ \\
\par
The remarks above and Theorem 4.2 show that Definition 4.1 fits with the non-Hom-associativity. The Hom-Akivis algebra constructed by Theorem 4.2 is said {\it associated} (with a given non-Hom-associative algebra). Starting from other considerations and using other notations, D. Yau has come to Theorem 4.2 above (see [21], Lemma 3.16).
\par
For the next results we need the following \\
\par
{\bf Definition 4.3.} Let $(\cal A , [-,-], [-,-,-], \alpha )$ and $(\tilde {\cal A} , \{-,- \}, \{-,-, - \}, \tilde {\alpha})$ be Hom-Akivis algebras. A morphism $\phi : \cal A \rightarrow \tilde {\cal A}$ of Hom-Akivis algebras is a linear map of $\mathbb K$-modules $\cal A$ and $\tilde {\cal A}$ such that
\par
$\phi ([x,y]) = \{ \phi (x), \phi (y) \}$,
\par
$\phi ([x,y,z]) = \{ \phi (x), \phi (y), \phi (z) \}$. \\
\par
For example, if take $(\cal A , [-,-], [-,-,-], \alpha )$ as a multiplicative Hom-Akivis algebra, then the twisting self-map $\alpha$ is itself an endomorphism of $(\cal A , [-,-], [-,-,-])$.
\par
The following result holds. \\
\par
{\bf Theorem 4.4.} {\it Let $(\cal A , [-,-], [-,-,-], \alpha )$ be a Hom-Akivis algebra and $\beta : \cal A \rightarrow \cal A$ a self-morphism of $(\cal A , [-,-], [-,-,-], \alpha )$. Define on $\cal A$ a bilinear operation $[-,-]_{\beta}$ and a trilinear operation $[-,-,-]_{\beta}$ by 
\par
$[x,y]_{\beta} := \beta ([x,y])$,
\par
$[x,y,z]_{\beta} := {\beta}^{2} ([x,y,z])$, \\
for all $x,y,z \in \cal A$, where $ {\beta}^{2} = \beta  \circ \beta $. Then ${\cal A}_{\beta} := (\cal A, [-,-]_{\beta} , [-,-,-]_{\beta}, \beta  \circ \alpha )$ is a Hom-Akivis algebra.
\par
Moreover, if $(\cal A , [-,-], [-,-,-], \alpha )$ is multiplcative and $\beta$ commutes with $\alpha$, then ${\cal A}_{\beta}$ is multiplicative}. \\
\par
{\sc Proof}: Clearly $[-,-]_{\beta}$ (resp. $[-,-,-]_{\beta}$) is a bilinear (resp. trilinear) map and the skew-symmetry of $[-,-]$ in $(\cal A , [-,-], [-,-,-], \alpha )$ implies the skew-symmetry of $[-,-]_{\beta}$ in ${\cal A}_{\beta}$.
\par
Next, we have (by the Hom-Akivis identity (4.1)), 
\par
$\sigma [[x,y]_{\beta} , (\beta \circ \alpha)(z) ]_{\beta} = {\beta}^{2}( \sigma [[x,y] , \alpha)(z)]) = {\beta}^{2}( \sigma [x,y,z] - \sigma [y,x,z] ) = \sigma ({\beta}^{2}([x,y,z]) - {\beta}^{2}([y,x,z])) = \sigma ([x,y,z]_{\beta} - [y,x,z]_{\beta})$ which means that ${\cal A}_{\beta}$ is a Hom-Akivis algebra.
\par
The second assertion is proved as follows:
\par
$[(\beta \circ \alpha)(x), (\beta \circ \alpha)(y)]_{\beta} = \beta ([(\beta \circ \alpha)(x), (\beta \circ \alpha)(y)]) =$ \\ $ \beta ([(\alpha \circ \beta)(x), (\alpha \circ \beta)(y)]) = (\beta \circ \alpha) ([\beta (x), \beta (y)]) = (\beta \circ \alpha) (\beta ([x,y])) = (\beta \circ \alpha) ([x,y]_{\beta})$ \\
and 
\par
$[ (\beta \circ \alpha)(x), (\beta \circ \alpha)(y), (\beta \circ \alpha)(z) ]_{\beta} = {\beta}^{2} ([ (\beta \circ \alpha)(x), (\beta \circ \alpha)(y), (\beta \circ \alpha)(z) ]) = {\beta}^{2} ([ (\alpha \circ \beta)(x), (\alpha \circ \beta)(y), (\alpha \circ \beta)(z) ]) = ({\beta}^{2} \circ \alpha) ([\beta (x), \beta (y), \beta (z)]) = $ \\ $(({\beta}^{2} \circ \alpha) \circ \beta ) ([x,y,z]) = ((\alpha \circ {\beta}^{2}) \circ \beta ) ([x,y,z]) = ((\alpha \circ \beta) \circ {\beta}^{2}) ([x,y,z]) =((\alpha \circ \beta)({\beta}^{2} ([x,y,z])) = (\beta \circ \alpha) ([x,y,z]_{\beta})$. \\
This completes the proof. \hfill $\square$ \\
\par
{\bf Corollary 4.5.} {\it If $(\cal A , [-,-], [-,-,-], \alpha )$ is a multiplicative Hom-Akivis algebra, then so is ${\cal A}_{\beta}$ }. \\
\par
{\sc Proof}: This is a special case of Theprem 4.4 when $\alpha = \beta$. \hfill $\square$ \\
\par
{\bf Corollary 4.6.} {\it Let $(\cal A , [-,-], [-,-,-])$ be an Akivis algebra and $\beta$ an endomorphism of $(\cal A , [-,-], [-,-,-])$. Define on $\cal A $ a bilinear operation $[-,-]_{\beta}$ and a trilinear operation $[-,-,-]_{\beta}$ by 
\par 
$[x,y]_{\beta} := [\beta (x), \beta (y)] \; \; (= \beta ([x,y]))$,
\par
$[x,y,z]_{\beta} := [{\beta}^{2} (x), {\beta}^{2} (y), {\beta}^{2} (z)] \; \; ( = {\beta}^{2} ([x,y,z]))$, \\
for all $x,y,z \in \cal A$, where ${\beta}^{2} = \beta \circ \beta$. Then $(\cal A , [-,-]_{\beta}, [-,-,-]_{\beta} , \beta)$ is a multiplicative Hom-Akivis algebra.
\par
Moreover, suppose that $\tilde {\cal A}$ is another Akivis algebra and that $\tilde \beta$ is an endomorphism of $\tilde {\cal A}$. If $f: {\cal A} \rightarrow \tilde {\cal A}$ is an Akivis algebra morphism satisfying $f \circ \beta = {\tilde \beta} \circ f$, then $f: (\cal A , [-,-]_{\beta}, [-,-,-]_{\beta} , \beta) \rightarrow (\tilde {\cal A} , [-,-]_{\tilde \beta}, [-,-,-]_{\tilde \beta} , \tilde \beta)$ is a morphism of multiplicative Hom-Akivis algebras}. \\
\par
{\sc Proof}: The first of this theorem is a special case of Theorem 4.4 above when $\alpha = id$. The second part is proved in the same way as in Theorem 4.4. For completeness, we repeat it as follows.
\par
$[f(x),f(y)]_{\tilde \beta} = {\tilde \beta}([f(x),f(y)]) = ({\tilde \beta} \circ f)([x,y]) = (f \circ \beta) ([x,y])$ \\ $= f([\beta (x), \beta (y)]) = f([x,y]_{\beta})$ \\
and
\par
$[f(x),f(y),f(z)]_{\tilde \beta} = {\tilde \beta}^{2}([f(x),f(y),f(z)]) = ({\tilde \beta}^{2} \circ f)([x,y,z])$ \\ $= ({\tilde \beta} \circ ({\tilde \beta} \circ f))([x,y,z]) = ({\tilde \beta} \circ (f \circ \beta))([x,y,z]) = (({\tilde \beta} \circ f) \circ \beta)([x,y,z])$ \\ $= ((f \circ \beta) \circ \beta)([x,y,z]) = (f \circ {\beta}^{2})([x,y,z]) = f({\beta}^{2}([x,y,z])) = f([x,y,z]_{\beta})$. \\
This completes the proof. \hfill $\square$ \\
\par
{\bf Remark.} (1) The gist of Theorem 4.4 is that the category of Hom-Akivis algebras is closed under twisting by self-morphisms. 
\par
(2) The Corollary 4.6 is the Akivis algebra analogue of a result in the Hom-Lie setting [19]. It shows how Hom-Akivis algebras can be constructed from Akivis algebras. This procedure was first given by Yau [19] in the construction of Hom-associative (resp. Hom-Lie) algebras starting from associative (resp. Lie) algebras. Such a procedure has been further extended to coalgebras [13] and to other systems (see, e.g., [5], [14]). \\
\\
{\bf 5. Hom-Malcev algebras from Hom-Akivis algebras} \\
\par
As for Akivis algebras, the notion of a Hom-Akivis algebra seems to be too wide in order to develop interesting specific results. For this purpose, it would be natural to consider some additional conditions and properties on Hom-Akivis algebras.
\par
In this section, we consider Hom-alternativity and Hom-flexibility in Hom-Akivis algebras. The main result here is that the Hom-Akivis algebra associated with a Hom-alternative algebra has a Hom-Malcev structure (this could be seen as another version of Theorem 3.8 in [21]).
\par
Since only the ternary operation of an Akivis algebra is involved in its alternativity or flexibility [3], we report these notions to Hom-Akivis algebras in the following \\
\par
{\bf Definition 5.1.} A Hom- Akivis algebra ${\cal A}_{\alpha} := (\cal A, [-,-], [-,-,-], \alpha )$ is said:
\par
(i) {\it Hom-flexible}, if $[x,y,x] = 0$, for all $x,y \in {\cal A}$;
\par
(ii) {\it Hom-alternative}, if $[-,-,-]$ is alternating (i.e. $[-,-,-]$ vanishes whenever any pair of variables are equal). \\
\par
{\bf Remark.} By linearization, as for associators in nonassociative algebras, one checks that 
\par
(1) the Hom-flexible law $[x,y,x] = 0$ in ${\cal A}_{\alpha}$ is equivalent to $ [x,y,z] = - [z,y,x]$;
\par
(2) the Hom-alternativity of ${\cal A}_{\alpha}$ is equivalent to its {\it left Hom-alternativity} ($[x,x,y] = 0$) and {\it right Hom-alternativity} $[y,x,x] = 0$ for all $x,y \in {\cal A}$. \\
\par
The following result is an immediate consequence of Theorem 4.2 and Definition 5.1. \\
\par
{\bf Theorem 5.2.} {\it Let $(A, \mu , \alpha )$ be a non-Hom-associative  algebra and $(A, [-,-], {\mathfrak as} (-,-,-) , \alpha )$ its associate Hom-Akivis algebra.
\par
(i) If $(A, \mu , \alpha )$ is Hom-flexible, then $(A, [-,-], {\mathfrak as} (-,-,-) , \alpha )$ is Hom-flexible.
\par
(ii) If $(A, \mu , \alpha )$ is Hom-alternative, then so is $(A, [-,-], {\mathfrak as} (-,-,-) , \alpha )$}. \hspace*{\fill} \hfill $\square$ \\
\par
We have the following characterization of Hom-Lie algebras in terms of Hom-Akivis algebras. \\
\par
{\bf Proposition 5.3.} {\it Let ${\cal A}_{\alpha} := (\cal A, [-,-], [-,-,-], \alpha )$ be a Hom-flexible Hom-Akivis algebra. Then ${\cal A}_{\alpha}$ is a Hom-lie algebra if, and only if $\sigma [x,y,z] = 0$, for all $x,y \in {\cal A}$}. \\
\par
{\sc Proof}: The Hom-Akivis identity (4.1) and the Hom-flexibility in ${\cal A}_{\alpha}$ imply 
\par
$\sigma [[x,y], \alpha (z) ] = 2 \; \sigma [x,y,z]$ \\
so that $\sigma [[x,y], \alpha (z) ] = 0$ if and only if $\sigma [x,y,z] = 0$ (recall that the ground field $\mathbb K$ is of characteristic 0). \hfill $\square$ \\
\par
The following result is a slight generalization of Proposition 3.17 in [21], which in turn generalizes a similar well-known result in alternative rings: \\
\par
{\bf Proposition 5.4.} {\it Let ${\cal A}_{\alpha} := (\cal A, [-,-], [-,-,-], \alpha )$ be a Hom-alternative Hom-Akivis algebra. Then }
\par
$\sigma [[x,y], \alpha (z) ] = 6 \; \sigma [x,y,z]$ \hfill (5.1) \\
{\it for all $x,y,z \in {\cal A}$}. \\
\par
{\sc Proof}: The application to (4.1) of the Hom-alternativity in ${\cal A}_{\alpha}$ gives the proof. \hfill $\square$ \\
\par
We now come to the main result of this section, which is Theorem 3.8 in [21] but from a point of view of Hom-Akivis algebras.
\par
In [21] D. Yau introduced the notion of a Hom-Malcev algebra: a {\it Hom-Malcev algebra} is a Hom-algebra $(A, [-,-], \alpha )$ such that the binary operation $[-,-]$ is skew-symmetric and that the identity \\
\par
$ {\Sigma}_{(\alpha (x), \alpha (y), z_x)}[[\alpha (x), \alpha (y)],\alpha (z_x)] = [\sigma [[x,y], \alpha (z)], {\alpha}^{2}(x)]$ \hfill (5.2) \\
\\
holds for all $x,y,z \in A$, where $ z_x := [x,z]$ and ${\Sigma}_{(\alpha (x), \alpha (y), z_x)}$ denotes the sum over cyclic permutation of $\alpha (x), \alpha (y)$, and $z_x$. The identity (5.2) is called the {\it Hom-Malcev identity}.
\par
Observe that when $\alpha = id $ then, by the skew-symmetry of $[-,-]$, the Hom-Malcev identity reduces to the Malcev identity ([15], [17]).
\par
The alternativity in Akivis algebras leads to Malcev algebras [3]. The Hom-version of this result is the following \\
\par
{\bf Theorem 5.5.} {\it Let $(A, \cdot , \alpha )$ be a Hom-alternative Hom-algebra and $(A, [-,-], {\mathfrak as}(-,-,-), \alpha)$ its associate Hom-Akivis algebra, where $[x,y] = x \cdot y - y \cdot x$ for all $x,y \in A$. Then $(A, [-,-], {\mathfrak as}(-,-,-), \alpha)$ is a Hom-Malcev algebra}. \\
\par
{\sc Proof}: From Theorem 4.2 we get that $(A, [-,-], {\mathfrak as}(-,-,-), \alpha)$ is Hom-alternative so that (5.1) implies \\
\par
$ {\Sigma}_{(\alpha (x), \alpha (y), z)}[[\alpha (x), \alpha (y)],\alpha (z)] = 6 \; {\mathfrak as}(\alpha (x), \alpha (y), z)$ \hfill (5.3) \\
\\
for all $x,y,z \in A$. Now, in (5.3) replace $z$ with $ z_x := [x,z]$ to get \\
\par
$ {\Sigma}_{(\alpha (x), \alpha (y), z_x)}[[\alpha (x), \alpha (y)],\alpha (z_x)] = 6 \; {\mathfrak as}(\alpha (x), \alpha (y), z_x)$ \hfill (5.4) \\
\\
But ${\mathfrak as}(\alpha (x), \alpha (y), z_x) = [{\mathfrak as}(x,y,z), {\alpha}^{2}(x)]$ in $(A, \cdot , \alpha )$ (see [21], Corollary 3.15) so that (5.4) reads \\
\par
$ {\Sigma}_{(\alpha (x), \alpha (y), z_x)}[[\alpha (x), \alpha (y)],\alpha (z_x)] = 6 \; [{\mathfrak as}(x,y,z), {\alpha}^{2}(x)]$ \\
\\
i.e., by (5.1) (viewing $[x,y,z]$ as ${\mathfrak as}(x,y,z)$), \\
\par
$ {\Sigma}_{(\alpha (x), \alpha (y), z_x)}[[\alpha (x), \alpha (y)],\alpha (z_x)] = [(\sigma [[x,y], \alpha (z)]), {\alpha}^{2}(x)]$ \\
\\
and one recognizes the Hom-Malcev identity (5.2). Therefore, we get that $(A, [-,-], {\mathfrak as}(-,-,-), \alpha)$ is a Hom-Malcev algebra. \hfill $\square$ \\
\par
{\bf Remark.} The procedure described in the proof of Theorem 5.5 somewhat repeats the one given by A.I. Maltsev in [15] when constructing Moufang-Lie algebras (now called Malcev algebras) from alternative algebras. \\
\par
{\bf Example 5.6.} By Theorem 4.2, the associate Hom-Akivis algebra of the Hom-flexible Hom-algebra of Example 3.4 is Hom-flexible (see Theorem 5.2). Moreover, Example 3.5 and Theorem 5.2 (see also Theorem 3.8 in [21]) imply that such a Hom-Akivis algebra is also a Hom-Malcev algebra. \\
\par
{\bf Acknowledgment.} I thank A. Makhlouf and D. Yau for their relevant comments and suggestions after reading the first draft of this paper, and for communicating with me about its specific topics.
\begin{center}
 \sc {References}
\end{center}
$[1]$ Akivis M.A., {\it Three-webs of multidimensional surfaces}, Trudy Geom. Sem. {\bf 2} (1969), 7-31 (Russian). \\
$[2]$ Akivis M.A., {\it Local differentiable quasigroups and three-webs of multidimensional surfaces}, in ``Studies in the Theory of Quasigroups and Loops'', 3-12, Izdat. ``Shtiintsa'', Kishinev, 1973 (Russian). \\
$[3]$ Akivis M.A., {\it Local algebras of a multidimensional three-web}, Siberian Math. J. {\bf 17} (1976), no. 1, 3-8. \\
$[4]$ Akivis M.A., Goldberg V.V. {\it Algebraic aspects of web geometry}, Comment. Math. Univ. Carolinae {\bf 41} (2000), no. 2, 205-236. \\
$[5]$ Ataguema H., Makhlouf A., Silvestrov S.D., {\it Generalization of n-ary Nambu algebra and beyond}, J. Math. Phys. {\bf 50} (2009), no. 8, 083501. \\
$[6]$ Fregier Y., Gohr A., {\it On Hom-type algebras}, arXiv:0903.3393v2 [math.RA] (2009). \\
$[7]$ Fregier Y., Gohr A., {\it On unitality conditions for Hom-associative algebras}, arXiv:0904.4874v2 [math.RA] (2009). \\
$[8]$ Gohr A., {\it On Hom-algebras with surjective twisting}, arXiv:0906.3270v3 [math.RA] (2009). \\
$[9]$ Hartwig J.T., Larsson D., Silvestrov S.D., {\it Deformation of Lie algebras using $\sigma$-derivations}, J. Algebra {\bf 295} (2006), 314-361. \\
$[10]$ Hofmann K.H., Strambach K., {\it Lie's fundamental theorems for local analytical loops}, Pacific J. Math. {\bf 123} (1986), no. 2, 301-327. \\
$[11]$ Kuzmin E. N., {\it Malcev algebras of dimension five over a field of zero characteristic}, Algebra i. Logika {\bf 9} (1970), 691-700. \\
$[12]$ Makhlouf A., Silvestrov S.D., {\it Hom-algebra structures}, J. Gen. Lie Theory Appl. {\bf 2} (2008), 51-64.\\
$[13]$ Makhlouf A., Silvestrov S.D., {\it Hom-algebras and Hom-coalgebras}, arXiv: 0811.0400v2 [math.RA] (2008). \\
$[14]$ Makhlouf A., {\it Hom-alternative algebras and Hom-Jordan algebras}, arXiv: 0909.0326v1 [math.RA] (2009). \\
$[15]$ Maltsev A.I., {\it Analytic loops}, Mat. Sb. {\bf 36} (1955), 569-576. \\
$[16]$ Myung H.C., {\it Malcev-Admissible Algebras}, Progress in Mathematics {\bf 64}, Birkh\"auser, Boston, MA, 1986. \\
$[17]$ Sagle A.A., {\it Malcev algebras}, Trans. Amer. Math. Soc. {\bf 101} (1961), 426-458. \\
$[18]$ Yau D., {\it Enveloping algebras of Hom-Lie algebras}, J. Gen. Lie Theory Appl. {\bf 2} (2008), 95-108. \\
$[19]$ Yau D., {\it Hom-algebras and homology}, J. Lie Theory {\bf 19} (2009), 409-421 (arXiv:0712.3515v2 [math.RA] (2008)). \\
$[20]$ Yau D., {\it Hom-Novikov algebras}, arXiv:0909.0726v1 [math.RA] (2009). \\
$[21]$ Yau D., {\it Hom-Maltsev, Hom-alternative, and Hom-Jordan algebras}, arXiv: 1002.3944v1 [math.RA] (2010).
\end{document}